\documentclass[12pt]{amsart}

\usepackage{amsmath,xspace,amssymb,mathrsfs}
\usepackage{color}

\input xy
\xyoption{all}
\xyoption{2cell}
\UseAllTwocells
\CompileMatrices

\newcommand{\Spec}{\operatorname{Spec}}
\renewcommand{\phi}{\varphi}

\newcommand{\Ker}{\operatorname{Ker}}
\newcommand{\Ima}{\operatorname{Im}}

\newcommand{\Tor}{\operatorname{Tor}}

\newcommand{\Max}{\operatorname{Max}}

\newtheorem{proposition}{Proposition}[section]
\newtheorem{lemma}[proposition]{Lemma} 

\newtheorem{theorem}[proposition]{Theorem}

\newtheorem{prop-def}[proposition]{Proposition and definition}

  \theoremstyle{definition}
  \newtheorem{definition}[proposition]{Definition}

\begin{document}

\title[Flat epimorphisms and pointwise rings]{On flat epimorphisms of rings and pointwise localizations}

\author[A. Tarizadeh]{Abolfazl Tarizadeh}
\address{ Department of Mathematics, Faculty of Basic Sciences,
University of Maragheh \\
P. O. Box 55136-553, Maragheh, Iran.
 }
\email{ebulfez1978@gmail.com}

\subjclass[2010]{14A05, 13B10, 13C11}
\keywords{Absolutely flat ring;
Flat epimorphism; Pointwise localiation.}

\begin{abstract} In this paper all rings are commutative. We prove some new results on flat epimorphisms of rings and pointwise localizations. Especially among them, it is proved that a ring $R$ is an absolutely flat (von-Neumann regular) ring if and only if it is isomorphic to the pointwise localization $R^{(-1)}R$, or equivalently, each $R-$algebra is $R-$flat. For a given minimal prime ideal $\mathfrak{p}$ of a ring $R$, the surjectivity of the canonical map $R\rightarrow R_{\mathfrak{p}}$ is characterized. Finally, we give a new proof to the fact that in a flat epimorphism of rings, the contraction-extension of an ideal equals the same ideal. 
\end{abstract}

\maketitle

\section{Introduction and Preliminaries}

In this paper, all rings are commutative. Let $\phi:R\rightarrow S$ be a morphism of rings and $J$ an ideal of $S$. Then it is easy to see that the contraction and then extension of $J$ under $\phi$ is contained in $J$, i.e., $J^{ce}\subseteq J$, see \cite[Proposition 1.17]{Atiyah}.
In this paper, in Theorem \ref{Theorem I new}, we give a new proof to the fact that if $\phi$ is a flat epimorphism of rings, then the equality holds. \\

Absolutely flat rings play a major role throughout this paper. We give two new characterizations for absolutely flat rings. The first one states that a given ring $R$ is absolutely flat if and only if each $R-$algebra is $R-$flat, see Theorem \ref{theorem 3}. The second characterization states that a ring $R$ is absolutely flat if and only if it is canonically isomorphic to the pointwise localization $R^{(-1)}R$, see Theorem \ref{lemma 234}. \\

By an epimorphism of rings $\phi:R\rightarrow S$ we mean it is an epimorphism in the category of commutative rings. It is important to notice that surjective ring maps are special cases of epimorphisms. As an example, the canonical ring map $\mathbb{Z}\rightarrow\mathbb{Q}$ is an epimorphism of rings which is not surjective. It is well known that a morphism of rings $R\rightarrow S$ is an epimorphism if and only if in the ring $S\otimes_{R}S$, $s\otimes1=1\otimes s$ for all $s\in S$. It is also well known that any faithfully flat epimorphism of rings is an isomorphism. In particular, an epimorphism of rings with source a field is an isomorphism if and only if the target is a nonzero ring. We refer the interested reader to \cite[Tag 04VM]{de Jong}, \cite{Lazard}, \cite{Olivier}, \cite{Roby}, \cite{Samuel}, \cite{ebulfez} and \cite{Tarizadeh} for a comprehensive discussion of epimorphisms of commutative rings. \\

By a flat epimorphism of rings we mean a ring map which is both a flat ring map and an epimorphism of rings. If $S$ is a multiplicative subset of a ring $R$, then the canonical ring map $R\rightarrow S^{-1}R$ is a typical example of flat epimorphisms of rings. \\

For a given ring $R$, the quotient ring $R/\mathfrak{N}$ is denoted by $R_{\mathrm{red}}$ where $\mathfrak{N}$ is the nil-radical of $R$. For any ring map $\phi:R\rightarrow S$ the induced map $R_{\mathrm{red}}\rightarrow S_{\mathrm{red}}$ is denoted by $\phi_{\mathrm{red}}$.\\

We shall freely use the above facts in this paper. \\

\section{Absolutely flat rings and epimorphisms of rings}

Recall that a ring $R$ is said to be an absolutely flat ring if each $R-$module is $R-$flat. It is well known that a ring $R$ is absolutely flat if and only if it is von-Neumann regular ring (i.e., each $r\in R$ can be written as $r=r^{2}s$ for some $s\in R$). In the following result we give a new and quite elementary proof for this well known fact.

\begin{theorem}\label{th 6 nice} Let $R$ be a ring. Then $R$ is an absolutely flat ring if and only if each $r\in R$ can be written as $r=r^{2}s$ for some $s\in R$.
\end{theorem}

{\bf Proof.} If $R$ is an absolutely flat ring then $R/I$ is $R-$flat where $I=(r)$. Then by \cite[Remark 2.2]{A. Tarizadeh acta}, $I=I^{2}$ and so $r=r^{2}s$ for some $s\in R$. To prove the reverse implication, it will be enough to show that for each $R-$module $M$ and for each ideal $I$ of $R$, then the canonical map $I\otimes_{R}M\rightarrow M$ which sends each pure tensor $a\otimes m$ of $I\otimes_{R}M$ into $am$ is injective. Assume $\sum\limits_{i=1}^{n}a_{i}m_{i}=0$ where $a_{i}\in I$ and $m_{i}\in M$ for all $i$. By hypothesis, each $a_{i}=r_{i}a^{2}_{i}$ for some $r_{i}\in R$. For $i=1,2$ we have $a_{i}=a_{i}b'$ where $b'=r_{1}a_{1}+r_{2}a_{2}-r_{1}r_{2}a_{1}a_{2}\in I$. Thus by induction on $n$, we may find some $b\in I$ such that $a_{i}=a_{i}b$ for all $i=1,\ldots,n$. It follows that
$\sum\limits_{i=1}^{n}a_{i}\otimes m_{i}=b\otimes(\sum\limits_{i=1}^{n}a_{i}m_{i})=0$. Hence, the above map is injective. $\Box$ \\

In the following result we provide a new characterization for absolutely flat rings. 

\begin{theorem}\label{theorem 3} Let $R$ be a ring. Then $R$ is absolutely flat if and only if each $R-$algebra is $R-$flat.
\end{theorem}

{\bf Proof.} The implication ``$\Rightarrow$" is clear. Conversely, let $M$ be a $R-$module. Then consider the ring $S=R\times M$ whose the addition and multiplication are defined as $(r,m)+(r',m')=(r+r',m+m')$ and $(r,m).(r',m')=(rr', rm'+r'm)$, respectively. Clearly $S$ is a commutative ring whose identity element is $(1,0)$ and the map $\phi: R\rightarrow S$ given by $r\rightsquigarrow(r,0)$ is a morphism of rings, (this construction is due to Nagata and in the literature, the ring $S$ is called the ``idealization'' or the trivial extension of $R$ by $M$). The $R-$module structure induced via $\phi$ on $S$ is the same as the usual $R-$module structure on the direct sum $R\oplus M$. By hypothesis, $\phi$ is a flat morphism. Thus $S=R\oplus M$ is a flat $R$-module.
It it well known that the direct sum of a family of $R$-modules is $R$-flat if and only if each factor is $R$-flat. Hence, $M$ is a flat $R-$module.  $\Box$ \\

\begin{lemma}\label{lemma revised 1} Let $S$ and $T$ be two multiplicative subsets of a ring $R$. Then $S^{-1}R\otimes_{R}T^{-1}R=0$ if and only if there exist $f\in S$ and $g\in T$ such that $fg=0$.
\end{lemma}

{\bf Proof.} It is proved exactly like \cite[Lemma 3.1]{Abolfazl 2020}. $\Box$ \\

\begin{theorem}\label{Theorem III 2020} Let $\mathfrak{p}$ be a minimal prime ideal of a ring $R$. Then the canonical map $\pi:R\rightarrow R_{\mathfrak{p}}$ is surjective if and only if the canonical map $R_{\mathfrak{m}}\rightarrow R_{\mathfrak{p}}$ is surjective for all $\mathfrak{m}\in\Max(R)\cap V(\mathfrak{p})$.
\end{theorem}

{\bf Proof.} The implication ``$\Rightarrow$" is clear, since the map $\pi$ factors as: $$\xymatrix{R\ar[r]&R_{\mathfrak{m}}\ar[r]&R_{\mathfrak{p}}.}$$
To see the converse, it suffices to show that the induced map $\pi_{\mathfrak{m}}:R_{\mathfrak{m}}\rightarrow
(R_{\mathfrak{p}})_{\mathfrak{m}}$ is surjective for all $\mathfrak{m}\in\Max(R)$. If $\mathfrak{p}\subseteq\mathfrak{m}$, then $(R_{\mathfrak{p}})_{\mathfrak{m}}\simeq R_{\mathfrak{p}}$. Thus by the hypothesis, $\pi_{\mathfrak{m}}$ is surjective. But if $\mathfrak{p}\nsubseteq\mathfrak{m}$, then choose $f\in\mathfrak{p}\setminus\mathfrak{m}$. Clearly $\mathfrak{p}R_{\mathfrak{p}}$ is the nil-radical of $R_{\mathfrak{p}}$. So there exists some $g\in R\setminus\mathfrak{p}$ such that $fg$ is nilpotent. Then by Lemma \ref{lemma revised 1} or by \cite[Lemma 3.1]{Abolfazl 2020}, $(R_{\mathfrak{p}})_{\mathfrak{m}}\simeq R_{\mathfrak{p}}\otimes_{R}R_{\mathfrak{m}}=0$. Hence, $\pi_{\mathfrak{m}}$ is surjective. $\Box$ \\

\begin{theorem}\label{Theorem I new} Let $\phi: R\rightarrow S$ be a flat epimorphism of rings. Then the following statements hold. \\
$\mathbf{(i)}$ If $\mathfrak{q}$ is a prime ideal of $S$, then the induced map $\phi_{\mathfrak{q}}: R_{\mathfrak{p}}\rightarrow S_{\mathfrak{q}}$ is an isomorphism of rings where $\mathfrak{p}=\phi^{-1}(\mathfrak{q})$.\\
$\mathbf{(ii)}$ If $J$ is an ideal of $S$, then $J^{ce}=J$. \\
$\mathbf{(iii)}$ The induced map $\phi^{\ast}:\Spec(S)\rightarrow\Spec(R)$ is a homeomorphism onto its image.\\
$\mathbf{(iv)}$ If $R$ is a Noetherian ring, then $S$ is as well.\\
$\mathbf{(v)}$ If $R$ is an Artinian ring, then $S$ is as well.\\
\end{theorem}

{\bf Proof.} $\mathbf{(i)}$: See \cite[Lemma 2.1]{A. Tarizadeh 5}. \\
$\mathbf{(ii)}$: Clearly $IS\subseteq J$ where $I:=\phi^{-1}(J)$. To see the reverse inclusion, we have $S/J\cong S/J\otimes_{R}S$ as $S-$modules. On the other hand, since $S$ is flat over $R$, thus from the exact sequence $\xymatrix{0 \ar[r] & R/I \ar[r]^{\overline{\phi}} & S/J}$ we obtain the following exact sequence $\xymatrix{0 \ar[r] & R/I\otimes_{R}S \ar[r]^{\overline{\phi}\otimes1_{S}} & S/J\otimes_{R}S}$.
Furthermore, the following diagram is commutative: $$\xymatrix{
R/I\otimes_{R}S \ar[r]^{\overline{\phi}\otimes1_{S}} \ar[d]^{} & S/J\otimes_{R}S \ar[d]^{} \\ S/IS\ar[r]^{} & S/J } $$
therefore $S/IS \rightarrow S/J$ is injective. Thus, $J\subseteq IS$. \\
$\mathbf{(iii)}$: By (ii), the function $\phi^{\ast}$ is a closed map onto its image.\\
$\mathbf{(iv)}$: Take an arbitrary ideal $J$ of $S$, then $I=\phi^{-1}(J)=(a_{1},...,a_{n})$ is a finitely generated ideal, since $R$ is a noetherian ring. By (ii), $J=IS=\big( \phi(a_{1}),...,\phi(a_{n})\big)$. Hence, $S$ is a Noetherian ring.\\
$\mathbf{(v)}$: Using (ii), then every descending chain of ideals of $S$ stabilizes. $\Box$  \\

Let $\phi: R\rightarrow S$ be a morphism of rings and let $J$ be the kernel of the canonical ring map $S\otimes_{R}S\rightarrow S$ given by $s\otimes s'\rightsquigarrow ss'$. Then it is well known that $J/J^{2}$ as $S-$module is canonically isomorphic to $\Omega_{R}(S)$, the module of K\"{a}hler differentials of $S$ over $R$, (it is also denoted by $\Omega_{S/R}$). In particular, $J$ is an idempotent ideal if and only if $\Omega_{R}(S)=0$. Then we provide a new proof to the following well known result.

\begin{theorem}\label{Theorem II new} A ring map $\phi: R\rightarrow S$ is an epimorphism of rings if and only if the following three conditions hold.\\
$\mathbf{(i)}$ The induced map $\phi^{\ast}:\Spec(S)\rightarrow\Spec(R)$ is injective.\\
$\mathbf{(ii)}$ If $\mathfrak{q}$ is a prime ideal of $S$, then the induced map $\kappa(\mathfrak{p})\rightarrow\kappa(\mathfrak{q})$ is an isomorphism where $\mathfrak{p}=\phi^{-1}(\mathfrak{q})$.\\
$\mathbf{(iii)}$ The kernel of the canonical ring map $p:S\otimes_{R}S\rightarrow S$ is a finitely generated and idempotent ideal.
\end{theorem}

{\bf Proof.} Assume $\phi:R\rightarrow S$ is an epimorphism of rings.
To prove (i), it will be enough to show that for each prime ideal $\mathfrak{p}$ of $R$, then $(\phi^{\ast})^{-1}(\mathfrak{p})$ has at most one element. It is well known that the fiber $(\phi^{\ast})^{-1}(\mathfrak{p})$ is homeomorphic to $\Spec\big(S\otimes_{R}\kappa(\mathfrak{p})\big)$. On the other hand, the base change $\kappa(\mathfrak{p})\rightarrow \kappa(\mathfrak{p})\otimes_{R}S$ is an epimorphism of rings. Hence, $\kappa(\mathfrak{p})\otimes_{R}S$ is a field whenever it is a non-zero ring. Therefore, $(\phi^{\ast})^{-1}(\mathfrak{p})$ has at most one element. To prove (ii), let $\mathfrak{q}$ be a prime ideal of $S$ laying over $\mathfrak{p}$, i.e., $\phi^{-1}(\mathfrak{q})=\mathfrak{p}$. Consider the following commutative diagram of rings: $$\xymatrix{
R \ar[r]^{\phi=epic} \ar[d]^{} & S \ar[d]^{epic} \\ \kappa(\mathfrak{p})\ar[r]^{} & \kappa(\mathfrak{q}).}$$
The composition  $\xymatrix{R \ar[r] & \kappa(\mathfrak{p}) \ar[r] & \kappa(\mathfrak{q})}$ is an epimorphism hence the induced ring map $\kappa(\mathfrak{p})\rightarrow\kappa(\mathfrak{q})$ is also an epimorphism (in fact it is an isomorphism). The statement (iii) is obvious, since $\phi$ is an epimorphism if and only if $p$ is an isomorphism. To prove the reverse implication we act as follows. Let $T$ be a reduced ring and let $f,g:S\rightarrow T$ be two ring maps such that $f\circ\phi=g\circ\phi$. We claim that $f=g$. We have $\phi^{\ast}\circ f^{\ast}=(f\circ\phi)^{\ast}=(g\circ\phi)^{\ast}=\phi^{\ast}\circ g^{\ast}$. Therefore $f^{\ast}=g^{\ast}$, since $\phi^{\ast}$ is injective. If $P$ is a prime ideal of $T$, then setting $\mathfrak{q}:=f^{-1}(\mathfrak{P})$, also setting $\mathfrak{p}:=\phi^{-1}(\mathfrak{q})$. Denote by $\widetilde{\phi}:\prod\limits_{P\in\Spec(T)}
\kappa(\mathfrak{p})\rightarrow\prod\limits_{P\in\Spec(T)}
\kappa(\mathfrak{q})$ the ring map induced via $\phi: R\rightarrow S$. By the hypothesis (ii), $\widetilde{\phi}$ is an isomorphism. Similarly above, denote by $\widetilde{f}, \widetilde{g}:S'=\prod\limits_{P\in\Spec(T)}
\kappa(\mathfrak{q})\rightarrow T'=\prod\limits_{P\in\Spec(T)}
\kappa(P)$ the ring maps induced by $f$ and $g$, respectively. From $f\circ\phi=g\circ\phi$ we conclude that $\widetilde{f}\circ\widetilde{\phi}=\widetilde{g}\circ\widetilde{\phi}$. Thus $\widetilde{f}=\widetilde{g}$, since $\widetilde{\phi}$ is an isomorphism. The following diagram is commutative: $$\xymatrix{
 S\ar[r]^{f,g} \ar[d]^{\rho} & T \ar[d]^{\rho'} \\ S'\ar[r]^{\widetilde{f},\widetilde{g}} & T' }$$
i.e. $\widetilde{f}\circ\rho=\rho'\circ f$ and $\widetilde{g}\circ\rho=\rho'\circ g$. Thus $\rho'\circ f=\rho'\circ g$. Since $T$ is reduced, therefore $\rho'$ is injective and so $f=g$. This establishes the claim. Now to conclude the assertion, since $i_{1}\circ\phi=i_{2}\circ\phi$ therefore $\eta\circ i_{1}\circ\phi=\eta\circ i_{2}\circ\phi$ where $i_{1}, i_{2}:S\rightarrow S\otimes_{R}S$ and
$\eta: S\otimes_{R}S\rightarrow(S\otimes_{R}S)_{\mathrm{red}}$ are the canonical maps. By applying what we have just proved for a morphism with target a reduced ring, we conclude that $\eta\circ i_{1}=\eta\circ i_{2}$. Therefore for each $s\in S$, $s\otimes1-1\otimes s$ is nilpotent. By the hypothesis (iii), $J=\Ker(p)$ is generated by a finite number of nilpotent elements of the form $s\otimes1-1\otimes s$. Hence,  $J$ is a nilpotent ideal. Thus $J=0$, since it is idempotent. Therefore, $\phi$ is an epimorphism of rings. $\Box$

\section{Pointwise rings with applications}

In this section we study the theory of pointwise localizations and some of its applications. This theory was originally introduced and studied during the S\'{e}minaire Samuel \cite{Samuel}. \\

If $R$ is an absolutely flat ring, then by Theorem \ref{th 6 nice}, each element $a\in R$ can be written as $a=a^{2}b$ for some $b\in R$. This leads us to the following definition.

\begin{definition} Let $R$ be a ring and let $a\in R$. If there is an element $b\in R$ such that $a=a^{2}b$ and $b=b^{2}a$, then $b$ is said to be a pointwise inverse of $a$.
\end{definition}

\begin{lemma}\label{lemma 9} Let $a,b\in R$. Then $b$ is a pointwise inverse of $a$ if and only if $a\in Ra^{2}$. Moreover, if $b$ is a pointwise inverse of $a$ then there is an idempotent element $e\in R$ such that $(e+a)(e+b)=1$. Finally, the pointwise inverse, if it exists, is unique.
\end{lemma}

{\bf Proof.} Suppose $a\in Ra^{2}$. We have $a=ra^{2}$ for some $r\in R$. Let $b=r^{2}a$. Then $b$ is a pointwise inverse of $a$. Clearly $e=1-ab$ is an idempotent element and $(e+a)(e+b)=1$.
Let $c\in R$ be another pointwise inverse of $a$. We have $b=ab^{2}=(ac)(ab^{2})=a^{2}c^{2}b=ac^{2}=c$. $\Box$ \\

The pointwise inverse of $a\in R$, if it exists, is usually denoted by $a^{(-1)}$. The pointwise inverse has appeared in the literature also
under other names, e.g. outer inverse or {2}-inverse. \\

\begin{lemma} Let $\phi: R\rightarrow S$ be a ring map. Suppose $a,b\in R$ have pointwise inverses in R. Then the pointwise inverses of $\phi(a)$ and $ab$ exist. Moreover $\phi(a)^{(-1)}=\phi(a^{(-1)})$ and $(ab)^{(-1)}=a^{(-1)}b^{(-1)}$.
\end{lemma}

{\bf Proof.} It is an easy exercise. $\Box$  \\

The following result establishes the universal property of the poinwise rings.

\begin{proposition}\label{pointwise} Let $R$ be a ring and let $S$ be a subset of $R$. Then there exist a ring $S^{(-1)}R$ and a canonical ring map $\eta: R\rightarrow S^{(-1)}R$ such that for each $s\in S$, the pointwise inverse of $\eta(s)$ in $S^{(-1)}R$ exists and the pair $(S^{(-1)}R, \eta)$ satisfies in the following universal property: if there is a ring map $\phi: R\rightarrow R'$ such that for each $s\in S$ the pointwise inverse of $\phi(s)$ in $R'$ exists then there is a unique ring map $\psi:S^{(-1)}R\rightarrow R'$ such that $\phi=\psi\circ\eta$.
\end{proposition}

{\bf Proof.} Consider the polynomial ring $A=R[x_{s} : s\in S]$ and let $S^{(-1)}R=A/I$ where the ideal $I$ is generated by elements of the form $sx_{s}^{2}-x_{s}$ and $s^{2}x_{s}-s$ with $s\in S$. Let $\eta: R\rightarrow S^{(-1)}R$ be the canonical ring map. For each $s\in S$, the element $x_{s}+I$ is the pointwise inverse of $\eta(s)=s+I$. Let $\phi:R\rightarrow R'$ be a ring map such that for each $s\in S$, the pointwise inverse of $\phi(s)$ exists in $R'$. By the universal property of the polynomial rings, there is a (unique) homomorphism of $R-$algebras $\widetilde{\phi}:R[x_{s} : s\in S]\rightarrow R'$ such that $x_{s}\rightsquigarrow\phi(s)^{(-1)}$ for all $s\in S$. We have $\widetilde{\phi}(I)=0$. Denote by $\psi:S^{(-1)}R\rightarrow R'$ the ring map induced by $\widetilde{\phi}$. Clearly $\psi$ is the unique ring homomorphism such that $\phi=\psi\circ\eta$. Because suppose there is another such ring map $\psi':S^{(-1)}R\rightarrow R'$. Then we have $\psi(x_{s}+I)=\widetilde{\phi}(x_{s})=\phi(s)^{(-1)}=
\psi'\big(\eta(s)\big)^{(-1)}=\psi'\big(\eta(s)^{(-1)}\big)=\psi'(x_{s}+I)$ for all $s\in S$. Therefore $\psi=\psi'$. $\Box$ \\

We call $S^{(-1)}R$ the pointwise localization of $R$ with respect to $S$. \\

\begin{proposition}\label{prop 6} Let $R$ be a ring and let $S$ be a subset of $R$. Then the following statements hold. \\
$\mathbf{(i)}$  The canonical ring map $\eta: R\rightarrow S^{(-1)}R$ is an epimorphism.\\
$\mathbf{(ii)}$ The map $\eta^{\ast}:\Spec\big(S^{(-1)}R\big)\rightarrow\Spec(R)$ is bijective.\\
$\mathbf{(iii)}$ For each $s\in S$, then $(\eta^{\ast})^{-1}\big(V(s)\big)$ is a clopen (both open and closed) subset of $\Spec\big(S^{(-1)}R\big)$.\\
$\mathbf{(iv)}$ The ring $S^{(-1)}R$ is nonzero if and only if $R$ is as well.\\
$\mathbf{(v)}$ $\Ker(\eta)\subseteq\mathfrak{N}$ where $\mathfrak{N}$ is the nil-radical of $R$.
\end{proposition}

{\bf Proof.} $\mathbf{(i)}:$ It follows from Proposition \ref{pointwise}.\\
$\mathbf{(ii)}:$ By Theorem \ref{Theorem II new}, the map $\eta^{\ast}$ is injective. To see surjectivity, let $\mathfrak{p}$ be a prime ideal of $R$ and consider the canonical ring map $\pi: R\rightarrow\kappa(\mathfrak{p})$. The image of every element of $R$ under $\pi$ has a pointwise inverse in $\kappa(\mathfrak{p})$. Thus, by Proposition \ref{pointwise}, there is a (unique) ring map $\psi:S^{(-1)}R\rightarrow\kappa(\mathfrak{p})$ such that $\pi=\psi\circ\eta$.
Then $\mathfrak{p}=\eta^{\ast}(\mathfrak{q})$ where
$\mathfrak{q}=\psi^{-1}(0)$.\\
$\mathbf{(iii)}:$ We have $(\eta^{\ast})^{-1}\big(V(s)\big)=V\big(\eta(s)\big)$. Moreover $V\big(\eta(s)\big)=D\big(1-\eta(s)\eta(s)^{(-1)}\big)$. \\
$\mathbf{(iv)}$ and $\mathbf{(v)}$: These are immediate consequences of (ii). $\Box$  \\

\begin{lemma}\label{lem 1} Let $\phi: R\rightarrow S$ be an epimorphism of rings where $S$ is a nonzero ring with trivial idempotents. Suppose $\phi(r)$ has a pointwise inverse in $S$ for all $r\in R$. Then $A:=\Ima(\phi)$ is an integral domain and $S$ is its field of fractions.
\end{lemma}

{\bf Proof.} Suppose $\phi(r)\phi(r')=0$ for some elements $r,r'\in R$. If $\phi(r)\neq0$ then $\phi(r)\phi(r)^{(-1)}=1$ since $\phi(r)\phi(r)^{(-1)}$ is an idempotent element. Therefore $A$ is an integral domain. Let $K$ be the field of fractions of $A$. Since every non-zero element of $A$ is invertible in $S$ therefore by the universal property of the localization, there is a (unique) ring map $\psi:K\rightarrow S$ such that $i=\psi\circ j$ where $i:A\rightarrow S$
and $j:A\rightarrow K$ are the canonical injections. The map $\phi$ factors as $\phi=i\circ\phi'$ where $\phi':R\rightarrow A$ is the ring map induced by $\phi$. Since $\phi$ is an epimorphism thus $i$ and so $\psi$ are epimorphisms. Hence, $\psi$ is an isomorphism. $\Box$ \\

\begin{theorem}\label{theorem 1} Let $R$ be a ring and let $\eta:R\rightarrow R'=R^{(-1)}R$ be the canonical ring map. Then the following statements hold. \\
$\mathbf{(i)}$ If $\mathfrak{q}$ is a prime ideal of $R'$, then $R'_{\mathfrak{q}}$ is canonically isomorphic to $\kappa(\mathfrak{p})$ where $\mathfrak{p}=\eta^{-1}(\mathfrak{q})$.\\
$\mathbf{(ii)}$ The ring $R^{(-1)}R$ is absolutely flat.
\end{theorem}

{\bf Proof.} $\mathbf{(i)}:$ For each prime ideal $\mathfrak{q}$ of $R^{(-1)}R$, the map: $$\xymatrix{R \ar[r]^{\eta} & R^{(-1)}R\ar[r]&R'_{\mathfrak{q}}}$$ satisfies all of the hypotheses of  Lemma \ref{lem 1}. Therefore $R'_{\mathfrak{q}}$ is a field. Now consider the following commutative diagram: $$\xymatrix{
R_{\mathfrak{p}} \ar[r]^{\eta_{\mathfrak{q}}=epic} \ar[d]^{}&R'_{\mathfrak{q}} \ar[d]^{\simeq} \\ \kappa(\mathfrak{p})\ar[r]^{} & \kappa(\mathfrak{q}) } $$
where $\mathfrak{p}=\eta^{\ast}(\mathfrak{q})$. The ring map $\kappa(\mathfrak{p})\rightarrow\kappa(\mathfrak{q})$ is an isomorphism, since it is an epimorphism.\\
$\mathbf{(ii)}:$ It is deduced from (i) and the fact that the absolutely
flatness is a local property. $\Box$  \\

By Proposition \ref{pointwise} and Theorem \ref{theorem 1}, the assignment $R\rightsquigarrow R^{(-1)}R$ is a covariant functor form the category of commutative rings into the category of absolutely flat rings. \\

\begin{lemma}\label{remark 1} Let $\phi: R\rightarrow S$ be a ring map, let $M$ and $N$ be $S-$modules and consider the canonical map $\eta:M\otimes_{R}N\rightarrow M\otimes_{S}N$ which maps each pure tensor $m\otimes_{R}n$ into $m\otimes_{S}n$. Then $\Ker(\eta)$ is generated by elements of the form $sm\otimes_{R}n-m\otimes_{R}sn$ with $s\in S\setminus\Ima(\phi)$, $m\in M$ and $n\in N$. In particular, if $\phi$ is an epimorphism of rings then $\eta$ is an isomorphism.
\end{lemma}

{\bf Proof.} Let $K$ be the $R-$submodule of $M\otimes_{R}N$ generated by elements of the form $sm\otimes_{R}n-m\otimes_{R}sn$ with $s\in S\setminus\Ima(\phi)$, $m\in M$ and $n\in N$. Clearly $K\subseteq\Ker(\eta)$. Consider the map
$\overline{\eta}:P=M\otimes_{R}N/K\rightarrow M\otimes_{S}N$ induced by $\eta$. We have $\Ker(\overline{\eta})=\Ker(\eta)/K$. The scalar multiplication $S\times P\rightarrow P$ which is defined on pure tensors by $s.(m\otimes_{R}n+K)=sm\otimes_{R}n+K$ is actually well-defined and puts a $S-$module structure over $P$. By the universal property of the tensor products, the $S-$bilinesr map $M\times N\rightarrow P$ defined by $(m,n)\rightsquigarrow m\otimes_{R}n+K$ induces a (unique) $S-$homomorphism $M\otimes_{S}N\rightarrow P$ which maps each pure tensor $m\otimes_{S}n$ into $m\otimes_{R}n+K$. This implies that $\overline{\eta}$ is bijective. Therefore $\Ker(\eta)=K$.  $\Box$ \\

\begin{lemma}\label{lem 2} Let $\phi: R\rightarrow S$ be a flat ring map which has a factorization $\xymatrix{ R \ar[r]^{\psi} & A \ar[r]^{\phi'} & S}$ such that $\phi'$ is an injective ring map and $\psi$ is an epimorphism of rings. Then $\phi'$ is a flat ring map.
\end{lemma}

{\bf Proof.}  For each $A-$module $M$, the canonical map $\eta_{M}:M\otimes_{R}S\rightarrow M\otimes_{A}S$ which maps each pure tensor $m\otimes_{R} s$ into $m\otimes_{A}s$ is injective because in $A\otimes_{R}A-$module $M\otimes_{R}S$ we have  $am\otimes_{R}s=(a\otimes_{R}1_{A}).(m\otimes_{R}s)=
(1_{A}\otimes_{R}a).(m\otimes_{R}s)=m\otimes_{R}a.s$ then apply Lemma \ref{remark 1}. In fact, it is bijective. Now suppose $\xymatrix{0 \ar[r] & N \ar[r]^{f} & M}$ is an exact sequence of $A-$modules. The following diagram is commutative:
$$\xymatrix{
 N\otimes_{R}S\ar[r]^{f\otimes_{R}1} \ar[d]^{\eta_{N}} & M\otimes_{R}S \ar[d]^{\eta_{M}} \\ N\otimes_{A}S \ar[r]^{f\otimes_{A}1} & M\otimes_{A}S} $$ and the map $f\otimes_{R}1$ is injective since $S$ is flat over $R$. Therefore $f\otimes_{A}1$ is injective as well. Hence, $S$ is a flat module over $A$, i.e., $\phi'$ is a flat ring map. $\Box$  \\

\begin{lemma}\label{simplified lemma} Let $\phi: R\rightarrow S$ be a flat epimorphism of rings. Then for each prime $\mathfrak{p}$ of $R$ we have either $\mathfrak{p}S=S$ or the canonical ring map $R_{\mathfrak{p}}\rightarrow T^{-1}S$ given by $r/s\rightsquigarrow\phi(r)/\phi(s)$ is an isomorphism where $T=\phi(R\setminus\mathfrak{p})$.
\end{lemma}

{\bf Proof.} Suppose $\mathfrak{p}S\neq S$ for some prime $\mathfrak{p}$. The canonical map $R_{\mathfrak{p}}\rightarrow T^{-1}S$ is a flat epimorphism because flat morphisms and epics are stable under base change and composition (recall that the ring $ T^{-1}S$ is canonically isomorphic to $S_{\mathfrak{p}}$). It is also faithfully flat since $\mathfrak{p}S\neq S$. Hence, it is an isomorphism. $\Box$ \\

It is worth mentioning that the converse of Lemma \ref{simplified lemma} also holds. \\

Theorems \ref{outstanding lemma} and \ref{coro 3} are well known, we provide new proofs for them.

\begin{theorem}\label{outstanding lemma} Let $\phi: R\rightarrow S$ be a flat epimorphism of rings. If $\phi_{\mathrm{red}}$ is surjective, then $\phi$ is as well.
\end{theorem}

{\bf Proof.} The map $\phi$ factors as $\xymatrix{R \ar[r]^{\pi\:\:\:\:\:\:\:\:\:\:\:} & R/\Ker(\phi)\ar[r]^{\:\:\:\:\:\:\:\:\:\:\:\:\phi'} & S}$
where $\pi$ is the canonical ring map and $\phi'$ is induced by $\phi$. We have $\Ima(\phi)=\Ima(\phi')$, $\phi'$ is an epimorphism and  $\phi'_{\mathrm{red}}$ is surjective. Moreover, by Lemma \ref{lem 2}, $\phi'$ is flat. Therefore, without loss of generality, we may assume that $\phi$ is injective. It follows that $\phi_{\mathrm{red}}$ is an isomorphism and so $\phi^{\ast}: \Spec(S)\rightarrow\Spec(R)$ is bijective. Therefore  $\mathfrak{p}S\neq S$ for all primes $\mathfrak{p}$ of $R$ and so by Lemma \ref{simplified lemma}, the canonical map $R_{\mathfrak{p}}\rightarrow S_{\mathfrak{p}}$ is bijective. It follows that $S/\phi(R)\otimes_{R}R_{\mathfrak{p}}=0$ for all primes $\mathfrak{p}$. $\Box$ \\

\begin{theorem}\label{coro 3} Let $\phi: R\rightarrow S$ be an epimorphism of rings such that $R$ is absolutely flat. Then $\phi$ is surjective.
\end{theorem}

{\bf Proof.} The map $\phi$ factors as $\xymatrix{R \ar[r]^{\pi\:\:\:\:\:\:\:\:\:\:\:} & R/\Ker(\phi)\ar[r]^{\:\:\:\:\:\:\:\:\:\:\:\:\phi'} & S}$
where $\pi$ is the canonical ring map and $\phi'$ is the injective ring map induced by $\phi$. The quotient ring $R/\Ker(\phi)$ is absolutely flat. Moreover, $\Ima(\phi)=\Ima(\phi')$ and yet $\phi'$ is an epimorphism. Hence, without loss of generality, we may assume that $\phi$ is injective. In this case, $\phi$ is a faithfully flat morphism. Because, suppose $S\otimes_{R}M=0$ for some $R-$module $M$. From the following short exact sequence of $R-$modules $$\xymatrix{0 \ar[r] & R\ar[r]^{\phi} & S\ar[r]^{\pi} & S/R \ar[r] & 0}$$ we obtain the following long exact sequence of $R-$modules
$\xymatrix{... \ar[r] &}$\\$\xymatrix{\Tor^{R}_{1}(S/R, M) \ar[r] &R\otimes_{R}M\ar[r]^{\phi\otimes1_{M}} & S\otimes_{R}M \ar[r]^{\pi\otimes1_{M}} & S/R\otimes_{R}M \ar[r] & 0}$.\\
But $\Tor^{R}_{1}(S/R,M)=0$ since $S/R$ is $R-$flat, see \cite[Theorem 7.2]{Rotman}. Thus $M\simeq R\otimes_{R}M=0$. Therefore $\phi$ is a faithfully flat epimorphism and so it is an isomorphism. This means that, in our factorization $\phi=\phi'\circ\pi$, $\phi'$ is an isomorphism. Therefore the original $\phi$ is surjective. $\Box$ \\

\begin{theorem}\label{lemma 234} Let $R$ be a ring. Then $R$ is an absolutely flat ring if and only if the canonical ring map $\eta:R\rightarrow R^{(-1)}R$ is an isomorphism.
\end{theorem}

{\bf Proof.} Suppose $R$ is absolutely flat. Then, by Theorem \ref{coro 3}, $\eta$ is surjective. Pick $a\in\Ker(\eta)$. By Theorem \ref{th 6 nice}, there exists some $b\in R$ such that $a=ba^{2}$. It follows that $a=b^{n-1}a^{n}$ for all $n\geqslant1$. But $a$ is a nilpotent element, see Proposition \ref{prop 6}. Therefore $a=0$. The reverse implication follows from Theorem \ref{theorem 1}. $\Box$ \\

\textbf{Acknowledgements.} The author would like to give sincere thanks to the referee for very careful reading of the paper. \\

\end{document}